\documentclass[leqno,12pt]{amsart}

\usepackage{amssymb,amscd}

\calclayout

\newtheorem{theorem}{Theorem}
\newtheorem{corollary}{Corollary}
\newtheorem{proposition}{Proposition}
\newtheorem{lemma}{Lemma}

\newcommand{\cL}{{\mathcal L}}
\newcommand{\cO}{{\mathcal O}}
\newcommand{\cX}{{\mathcal X}}

\newcommand{\mA}{{\mathbb A}}
\newcommand{\mC}{{\mathbb C}}
\newcommand{\mG}{{\mathbb G}}
\newcommand{\mP}{{\mathbb P}}
\newcommand{\mQ}{{\mathbb Q}}
\newcommand{\mR}{{\mathbb R}}
\newcommand{\mZ}{{\mathbb Z}}

\newcommand{\diag}{\operatorname{diag}}
\renewcommand{\div}{\operatorname{div}}
\newcommand{\reg}{\operatorname{reg}}
\newcommand{\Pic}{\operatorname{Pic}}
\newcommand{\Proj}{\operatorname{Proj}}
\newcommand{\PGL}{\operatorname{PGL}}
\newcommand{\SL}{\operatorname{SL}}

\title{The cone of effective one--cycles of certain $G$--varieties}

\author{Michel~Brion}
\address{Universit\'e de Grenoble I\\
D\'epartement de Math\'ematiques\\
Institut Fourier, UMR 5582 du CNRS\\
38402 Saint-Martin d'H\`eres Cedex, France}
\email{Michel.Brion@ujf-grenoble.fr}

\begin{document}

\begin{abstract}
Let $X$ be a normal projective variety admitting an action of a
semisimple group with a unique closed orbit. We construct finitely
many rational curves in $X$, all having a common point, such that
every effective one--cycle on $X$ is rationally equivalent to a
unique linear combination of these curves with non--negative rational 
coefficients. When $X$ is nonsingular, these curves are projective
lines, and they generate the integral Chow group of one--cycles.
\end{abstract}

\maketitle

\medskip

\hfill{\it To C.~S.~Seshadri for his 70th birthday}

\section*{Introduction}

One associates to any complete algebraic variety $X$ the group $N_1(X)$
of one--cycles on $X$ modulo numerical equivalence; this is a free
abelian group of finite rank \cite[Example 19.1.4]{Fu98}. In the
corresponding real vector space $N_1(X)_{\mR}$, the convex cone
generated by the classes of (closed irreducible) curves in $X$ is the
{\it cone of effective one--cycles}, denoted by $NE(X)$.

\smallskip

Dually, one has the group $N^1(X)$ of Cartier divisors on $X$ modulo
numerical equivalence, endowed with a non--degenerate pairing 
$N^1(X)\times N_1(X)\rightarrow\mZ$ via intersection numbers of Cartier
divisors with curves. Again, $N^1(X)$ is a free abelian group of
finite rank; in the corresponding real vector space $N^1(X)_{\mR}$,
the dual cone of $NE(X)$ is the {\it cone of numerically effective
divisors} (called {\it nef} for brevity). The interior of the nef cone
is the ample cone, by Seshadri's criterion \cite[Theorem I.7.1]{Ha70}. 

\smallskip

The cone $NE(X)$ encodes much information about morphisms with source
$X$. Specifically, any surjective separable morphism 
$f:X\rightarrow Y$ with connected fibers, where $Y$ is a normal
projective variety, is uniquely determined by the subcone
$NE(f)\subseteq NE(X)$ generated by the classes of curves lying in
fibers of $f$. Moreover, $NE(f)$ is a face of the convex cone $NE(X)$. 

\smallskip

However, the structure of $NE(X)$ may be quite complicated: this cone
may not be closed, its closure may not be polyhedral, and certain
faces may not arise from morphisms with source $X$, due to the
existence of nef divisors having no globally generated positive
multiple. Only the part of $NE(X)$ where the canonical divisor is 
negative is well--understood, if $X$ is a nonsingular projective
variety in characteristic zero (see e.g.~\cite{CKM88}).

\smallskip

In this note, we consider certain varieties where the structure
of the whole cone of effective one--cycles turns out to be very
simple. These are the normal projective varieties where a semisimple
group acts with a unique closed orbit. For such a variety $X$, we show
that the cone $NE(X)$ is generated by the closures of positive strata
of dimension one, for an appropriate Bialynicki--Birula decomposition
of $X$ (recalled in Section 1). These closures are rational curves
passing through a common point, the sink of the decomposition; their
classes form a basis of the rational vector space $N_1(X)_{\mQ}$, and
the latter is isomorphic to the rational Chow group of one--cycles. 
If, in addition, $X$ is nonsingular, then so are our rational curves,
and their classes form a basis of the group $N_1(X)$ (Theorem
\ref{curves} and Corollary \ref{equivalences}). 

\smallskip

Moreover, the cone of nef divisors is generated by the closures of
negative strata of codimension one. These are globally generated
Cartier divisors, and their classes in $N^1(X)_{\mQ}$ form the dual
basis to our basis of $N_1(X)_{\mQ}$; finally, $N^1(X)$ is
isomorphic to the Picard group. As a consequence, every nef divisor on
$X$ is globally generated (Theorem \ref{divisors}).

\smallskip

The simplest examples of projective varieties where a semisimple group
acts with a unique closed orbit are of course flag varieties. For
these, our results are well--known: the cone of effective one--cycles
is freely generated by the classes of Schubert varieties of dimension
one, while the classes of opposite Schubert varieties of codimension
one form the dual basis of the nef cone. Moreover, numerical and
rational equivalence coincide. 

\smallskip

These results were generalized in \cite{Br93} and \cite{Br97} to all
projective simple spherical varieties (i.e., normal projective
varieties where a connected reductive group acts with a unique closed
orbit, and where a Borel subgroup acts with a dense orbit), in
characteristic zero. There the main tool was the classification of
embeddings of spherical homogeneous spaces. In the present note, the
use of the Bialynicki--Birula decomposition (for possibly singular
varieties) yields more general results.

\smallskip

Another interesting class of examples consists of orbit closures of a
finite set of ordered points of the projective line, under the
diagonal action of $\PGL(2)$. Clearly, these orbit closures are
projective varieties where $\PGL(2)$ acts with a unique closed orbit;
and they turn out to be normal, as shown by Iozzi and Poritz (in
characteristic zero; see \cite{IP99}). In the final section of this
note, we show that their cone of nef divisors is generated by the
pull--backs of points under the various projections to the projective
line, while the cone of effective curves is generated by the
projective lines mapped to points under all but one projection
(Proposition \ref{PGL}). We also obtain another proof of the normality
of these orbit closures, and we show that they yield all simple
complete embeddings of $\PGL(2)$, in arbitrary characteristics.

\smallskip

Note finally that the duality between closures of positive strata of
dimension one and of negative strata of codimension one, is very
seldom satisfied by (say) nonsingular projective varieties where
a torus acts with only finitely many fixed points. One checks, for
example, that this duality fails for all toric surfaces having at
least $5$ fixed points.

\smallskip

\noindent
{\it Acknowledgments.} Many thanks to L.~Bonavero and S.~Druel for
very hepful discussions and comments; their joint work \cite{BCDD02}
(with C.~Casagrande and O.~Debarre) on a conjecture of Mukai has been
the starting point of this note.

\section{The cone of numerically effective divisors}

We begin by fixing notation and recalling some general results on
group actions, to be used throughout this work. 

All algebraic varieties and groups are defined over an algebraically
closed field $K$ of arbitrary characteristic. 

Let $G$ be a connected reductive algebraic group. Let 
$B, B^-\subseteq G$ be opposite Borel subgroups with common torus $T$
and unipotent radicals $U=R_u(B), U^-=R_u(B^-)$; then $B=TU$ and
$B^-=TU^-$. Let $\cX^*(T)$ be the character group of $T$ (the weight
lattice of $G$), and $\cX_*(T)$ the group of one--parameter subgroups
of $T$ (the dual lattice of $\cX^*(T)$).

A $G$--{\it module} is a the space of a rational, finite--dimensional
$G$--representation. Important examples of $G$--modules are the
{\it induced modules} $H^0(\omega)$, where $\omega\in\cX^*(T)$: these
consist of all regular functions $f$ on $G$ such that
$f(gtu)=\omega(t^{-1})f(g)$ for all $g\in G$, $t\in T$ and $u\in U$. 
The $G$--action is defined by $(g\cdot f)(h)=f(g^{-1}h)$. Recall that
$H^0(\omega)$ is non--zero if and only if $\omega$ is dominant; 
then the set of $U$--fixed points in $H^0(\omega)$ is a line where
$T$ acts with weight $\omega$ (see \cite[II.2]{Ja87} for these
results, and for more on induced modules).

A {\it variety} is a separated integral scheme of finite type over
$K$. For a connected linear algebraic group $H$, a $H$--{\it variety} 
is a variety $X$ endowed with an algebraic action of $H$. The 
fixed point subset of $H$ in $X$ is denoted by $X^H$ and regarded as a
reduced closed subscheme of $X$. An example of a $H$--variety is the
projectivization $\mP(M)$ of a $H$--module $M$; an $H$--variety $X$ is
$H$--{\it quasi--projective} if it admits an equivariant embedding
into some $\mP(M)$. Equivalently, $X$ admits an ample
$H$--linearizable invertible sheaf. Any normal $H$--variety $X$ admits
a covering by $H$--quasi--projective open subsets; as a consequence,
$X$ is $H$--quasi--projective if it is quasi--projective, or if it
contains a unique closed $H$--orbit \cite[Theorem 1]{KKLV89}.

Let $\cL$ be an invertible sheaf on a normal $H$--variety. If the
Picard group of $H$ is trivial, then $\cL$ admits a
$H$--linearization; in the general case, there exists a finite
covering $\tilde H\rightarrow H$ of algebraic groups such that 
$\tilde H$ has a trivial Picard group. Thus, $\cL$ is 
$\tilde H$--linearizable, and some positive power $\cL^n$ is
$H$--linearizable (for all this, see \cite[I.3]{MFK94} and also
\cite{KKLV89}, \cite{KKV89}).

\medskip

Next we obtain a criterion for a $G$--projective variety to contain a
unique closed orbit.

\begin{lemma}\label{point}
The following conditions are equivalent for a $G$--projective variety
$X$:

\begin{enumerate}

\item
$X$ contains a unique closed $G$--orbit.

\item
$X^U$ consists of a unique point.

\item
There exists a finite equivariant morphism
$f:X\rightarrow\mP(H^0(\omega))$, for some dominant weight $\omega$.

\end{enumerate}

Then $X$ is fixed pointwise by the connected centre $Z(G)^0$, so
that $G$ acts on $X$ via its semisimple quotient $G/Z(G)^0$.
\end{lemma}

\begin{proof}
(1)$\Rightarrow$(2) By assumption, $X^B$ consists of a unique point
$x$. In other words, $T$ acts on $X^U$ with $x$ as its unique fixed
point. Since $X^U$ is a $T$--projective variety, it follows that
$X^U=\{x\}$.

(2)$\Rightarrow$(1) follows from Borel's fixed point theorem.

(1)$\Rightarrow$(3) Let $\cL$ be a $G$--linearized, very ample
invertible sheaf on $X$; let $Y\subseteq X$ be the closed
$G$--orbit. Then the $G$--module $H^0(Y,\cL)$ contains a
$B$--eigenvector $\sigma$. By \cite[Korollar 2.3]{Kn93}, there exist a 
positive integer $n$ such that $\sigma^n\in H^0(Y,\cL^n)$
lifts to a $B$--eigenvector $\tau\in H^0(X,\cL^n)$. Let
$N\subseteq H^0(X,\cL^n)$ be the $G$--submodule generated
by $\tau$. Then the set of common zeroes in $X$ of all sections in $N$
is closed, $G$--stable, and does not contain $Y$. Hence this set is
empty, since $Y$ is the unique closed orbit. So the canonical rational
map $X -\rightarrow \mP(N^*)$ is in fact a finite equivariant
morphism. Moreover, the $G$--module $N$ being generated by a
$B$--eigenvector, it is a quotient of a universal highest weight
module; thus, its dual $N^*$ is a submodule of some $H^0(\omega)$, by
\cite[II.2.13]{Ja87}.  

(3)$\Rightarrow$(2) Since $U$ is unipotent connected and
$H^0(\omega)^U$ is a line, $\mP(H^0(\omega))^U$ consists of a
unique point; thus, the set $X^U$ is finite. On the other hand, $X^U$
is connected by \cite{Ho69}; hence it consists of a unique point.

Finally, the centre $Z(G)$ acts on $H^0(\omega)$ by scalars
\cite[II.2.8]{Ja87}, and hence fixes $\mP(H^0(\omega))$ pointwise. 
Since $f$ is finite, it follows that $Z(G)^0$ fixes $X$ pointwise. 
\end{proof}

\medskip

\noindent
{\bf Definition.} A $G$--projective variety satisfying one of the
conditions of Lemma \ref{point} is called {\it simple}.

\medskip

We will construct generators of the nef cone of a simple $G$--projective 
variety, by using the {\it Bialynicki--Birula decomposition} that we
now review, after \cite{Bi73} (in the nonsingular case) and
\cite{Ko78} (in the general case).

Let $X$ be a $T$--projective variety and let $\lambda\in\cX_*(T)$. The
multiplicative group $\mG_m$ acts on $X$ via 
$\lambda:\mG_m\rightarrow T$, and we denote by $X^{\lambda}$ the
corresponding fixed point set. 
Given $x\in X$, the map $\mG_m\rightarrow X,~t\mapsto\lambda(t)\cdot x$
extends uniquely to a morphism $\mP^1\rightarrow X$; this defines 
$\lim_{t\to 0}\lambda(t)\cdot x$ and 
$\lim_{t\to\infty}\lambda(t)\cdot x$. Both are $\lambda$--fixed
points, distinct unless $x\in X^{\lambda}$. 

Recall also that $X^{\lambda}$ equals $X^T$ for any $\lambda\in X_*(T)$
outside a finite union of proper subgroups; such a $\lambda$ is
called {\it regular}. 

For any $\lambda\in \cX_*(T)$ and any closed subset 
$Y\subseteq X^{\lambda}$, let 
$$
X^+(Y)=\{x\in X~\vert~ \lim_{t\to 0} \lambda(t)\cdot x\in Y\}.
$$
Then $X^+(Y)$ is a locally closed $T$--invariant subset of $X$, and
the map
$$
p^+:X^+(Y)\to Y,~x\mapsto \lim_{t\to 0} \lambda(t)\cdot x
$$
is a surjective affine $T$--invariant morphism (these facts are easily
checked in the ambient projectivization of a $T$-module). Moreover,
$X$ is the disjoint union of the subsets $X^+(Y)$, where $Y$ runs over
all connected components of $X^T$. These subsets are called the 
{\it positive strata}. 

As a consequence, there exists a unique open positive stratum $X^+(Y)$.
Then $Y=Y^+$ is irreducible; it is called the {\it source} of $X$ for
the $\lambda$--action. 

Likewise, $X$ is the disjoint union of the 
$$
X^-(Y)=\{x\in X~\vert~ \lim_{t\to \infty} \lambda(t)\cdot x\in Y\},
$$
where $Y$ runs over the connected components of $X^{\lambda}$. The
unique component $Y=Y^-$ such that $X^-(Y)$ is open in $X$ is called
the {\it sink}.

Note that $X^+(Y)\cap X^-(Y)=Y$ for any $Y\subseteq X^{\lambda}$, so
that $X^-(x)=\{x\}$ for all $x\in Y^+$. Moreover, by \cite{Ko78},
one has for any $x\in X^{\lambda}$: 
$$
\dim X-\dim_x X^{\lambda}\le \dim_x X^+(x) + \dim_x X^-(x). \leqno(*)
$$
As a consequence, the source of $X$
consists of those $x\in X^{\lambda}$ such that $X^-(x)=\{x\}$
(indeed, if $X^-(x)=\{x\}$, then $X^-(x')=(p^-)^{-1}(x')$ equals
$\{x'\}$ for all $x'$ in a neighborhood of $x$ in $X^{\lambda}$, since
$\{x'\}$ is the unique closed $\mG_m$--orbit in its fiber under $p^-$.
Now ($*$) implies that $\dim X = \dim X^+(Y)$ for any irreducible
component $Y\subseteq X^{\lambda}$ containing $x$.)

If, in addition, $X$ is nonsingular, then so are $X^{\lambda}$ and
$X^+(Y)$, $X^-(Y)$ for all components $Y$ of $X^{\lambda}$. Moreover,
$X^+(Y)$ and $X^-(Y)$ intersect transversally along $Y$, and the
morphisms $p^+$, $p^-$ are affine bundles, by \cite{Bi73}.

Next we recall the definition of certain subgroups of $G$ associated
with a given $\lambda\in\cX_*(T)$; see \cite[II.2.6]{MFK94} for details. 
Let
$$
G(\lambda)=\{g\in G~\vert~ \lambda(t)g\lambda(t^{-1}) 
\text{ has a limit in $G$ as }t\to 0 \},
$$
so that 
$$
G(-\lambda)=\{g\in G~\vert~ \lambda(t)g\lambda(t^{-1}) 
\text{ has a limit in $G$ as }t\to \infty \}.
$$
Then $G(\lambda)$ and $G(-\lambda)$ are opposite
parabolic subgroups of $G$, with common Levi subgroup the centralizer
of $\lambda$. Moreover, the unipotent radical of $G(\lambda)$ is 
$$
R_uG(\lambda)=\{g\in G~\vert~ 
\lim_{t\to 0}\lambda(t) g\lambda(t^{-1}) =1\}.
$$
Finally, $G(\lambda)$ equals $B$ if and only if $\lambda$ lies in the
interior of the Weyl chamber associated with $B$; then
$G(-\lambda)=B^-$.

Consider now a $G$--projective variety $X$. Then we may choose a
regular $\lambda\in\cX_*(T)$ such that $G(\lambda)=B$. Hence $X^+(Y)$
is $B$--invariant for any closed subset $Y\subseteq X^{\lambda} = X^T$ 
(since 
$\lambda(t)g\cdot x = \lambda(t)g\lambda(t^{-1})\lambda(t)\cdot x$). 
Likewise, $X^-(Y)$ is $B^-$--invariant. This implies at once that the
sink $Y^-$ is fixed pointwise by $B$. In particular, we obtain

\begin{lemma}\label{sink}
Let $X$ be a simple $G$--projective variety. Then the sink consists of
the unique $B$--fixed point $x^-$, and the source consists of the
unique $B^-$--fixed point $x^+$.
\end{lemma}

Put for simplicity $X^-=X^-(x^-)$. Then $X^-$ is an open affine
$B^-$--invariant neighborhood of $x^-$; one easily checks that such a
neighborhood is unique. Moreover, $\mG_m$ acts on the algebra of
regular functions $R= H^0(X^-,\cO_{X^-})$ via 
$(t\cdot f)(x)=f(\lambda(t^{-1})x)$, and this yields a positive
grading of $R$. By Nakayama's Lemma, any graded invertible $R$--module
is generated by a homogeneous element, unique up to scalar. This
implies

\begin{lemma}\label{principal} 
With the preceding notation, the group of isomorphism classes of
$T$--linearized invertible sheaves on $X^-$ is isomorphic to
$\cX^*(T)$, via pull--back to the fixed point $x^-$.
\end{lemma}

We will also need the following variant of a result of Knop
\cite[Lemma 2.2]{Kn94}.

\begin{lemma}\label{generated}
Let $H$ be a connected linear algebraic group, $X$ a normal
$H$--variety, and $D$ an effective Weil divisor on $X$ whose support
contains no $H$--orbit. Then $D$ is a globally generated Cartier
divisor.
\end{lemma}

\begin{proof}
Let $i:X^{\reg}\rightarrow X$ be the inclusion of the nonsingular
locus. Then the sheaf $i^*\cO_X(D)$ is invertible. Replacing $H$ by a
finite cover, we may assume that this sheaf is $H$--linearized. Since
$\cO_X(D)\cong i_*i^*\cO_X(D)$( by normality of $X$), it follows that
$\cO_X(D)$ is $H$--linearized as well.

Let $Y\subseteq X$ be the subset of all points where $\cO_X(D)$ is not
invertible. Then $Y$ is closed, $H$--invariant, and contained in the
support of $D$; thus, $Y$ contains no orbit. Hence $Y$ is empty, and
$D$ is Cartier. Likewise, the base locus of $D$ is empty, so that $D$
is globally generated.
\end{proof}

We now come to our first main result.

\begin{theorem}\label{divisors}
Let $X$ be a normal, simple $G$--projective variety. Let $x^-\in X$ be
the $B$--fixed point, $X^-\subseteq X$ its unique $B^-$--invariant open
affine neighborhood, and $D_1,\ldots, D_r$ the irreducible components
of $X \setminus X^-$. Then the following hold:

\begin{enumerate}

\item 
$D_1,\ldots,D_r$ are globally generated Cartier divisors. Their
linear equivalence classes form a basis of the Picard group of $X$. 

\item
Every ample divisor on $X$ is linearly equivalent to a
unique linear combination of $D_1,\ldots,D_r$ with positive
integer coefficients. 

\item
Every nef divisor on $X$ is linearly equivalent to a
unique linear combination of $D_1,\ldots,D_r$ with non--negative
integer coefficients. 

\end{enumerate}

\end{theorem}

\begin{proof}
(1) The first assertion follows from Lemma \ref{generated} (applied
to $H=B$). For the second assertion, let $D$ be a Cartier divisor on
$X$. Then the invertible sheaf $\cO_X(D)$ admits a
$T$--linearization; by Lemma \ref{principal}, it follows that
the pull--back of $\cO_X(D)$ to $X^-$ is trivial (as an invertible
sheaf). Thus, $D$ is linearly equivalent to a linear combination of 
$D_1,\ldots, D_r$ with integer coefficients. These coefficients are
unique, since every regular invertible function on $X^-$ is constant.

(2) Let $D$ be an ample divisor on $X$. There exists a
positive integer $n$ such that the invertible sheaf $\cO_X(nD)$ is
very ample and $G$--linearized. Then the $G$--module 
$H^0(G\cdot x^-,\cO_X(nD))$ contains a $B^-$--eigenvector
$\sigma$. Note that $\sigma(x^-)\ne 0$, since $\sigma\ne 0$ and
$B^-\cdot x^- = B^- B\cdot x^-$ is dense in $G\cdot x^-$. Replacing
$n$ by a positive multiple, we may also assume that $\sigma$ lifts to
$\tau\in H^0(X,\cO_X(nD))$; then we may further assume that $\tau$
is a $T$--eigenvector. Since $\tau(x^-)\ne 0$, it follows that $\tau$
has no zero on $X^-$, so that $nD$ is linearly equivalent to a linear
combination of $D_1,\ldots,D_r$ with non--negative integer coefficients.

Since the divisor $mD-D_1-\cdots -D_r$ is ample for large $m$, this
shows that $D$ is linearly equivalent to a linear combination of
$D_1,\ldots,D_r$ with positive rational coefficients. Together with
(1), this implies our statement.

(3) Let $E$ be a nef divisor on $X$. Since $D_1+\cdots+D_r$ is
ample (by (2)), the divisor $E+D_1+\cdots+D_r$ is ample as
well. Applying again (2) completes the proof.
\end{proof}

In particular, the Picard group $\Pic(X)$ is a free abelian group
of finite rank. Together with \cite[Examples 19.1.2, 19.3.3]{Fu98},
this implies 

\begin{corollary}\label{picard}
For any normal, simple $G$--projective variety $X$, the natural map
$\Pic(X)\rightarrow N^1(X)$ is an isomorphism. In other words,
rational and numerical equivalence coincide for Cartier divisors.
\end{corollary}

\section{The cone of effective one--cycles}

We will construct generators of the cone of effective one--cycles of a
simple $G$--projective variety $X$, that are dual to the generators of
the nef cone obtained in Theorem \ref{divisors}. For this, recall that
every effective cycle on $X$ is rationally equivalent to a linear
combination of $B$--invariant subvarieties with positive coefficients,
by \cite{Hi84}. Thus, we will study the $B$--invariant curves in
$X$; we begin with the following easy result, generalizing 
\cite[Proposition 2.1]{Br93} to arbitrary characteristics.

\begin{lemma}\label{rational}
Let $C$ be a $B$--invariant curve in a simple $G$--projective
variety. Then $C$ is fixed pointwise by the radical of a unique
minimal parabolic subgroup $P(C)\supset B$, and 
$C=\overline{B\cdot x}$ for a unique $x\in C^T$. Moreover, the 
normalization of $C$ is isomorphic to $\mP^1$, and the normalization
map is bijective.
\end{lemma}

\begin{proof}
By assumption, $C^B$ consists of a unique point. And since $C$ is
$T$--invariant, $C\setminus C^B$ contains a unique $T$--fixed point
$x$. The reduced isotropy group $B_x$ equals $U_x T$, where $U_x$ is a
closed reduced subgroup of codimension one in $U$, normalized by
$T$. It follows easily that $U_x=R_u(P)$ for a unique minimal
parabolic subgroup $P\supset B$. Then the radical $R(P)$ is a normal
subgroup of $B$ contained in $B_x$. Thus, $R(P)$ fixes pointwise 
$\overline{B\cdot x}=C$, and $R_u(P)$ is the kernel of the $U$--action 
in $C$.

Since $C$ consists of a dense $B$--orbit and a unique fixed point, it
is rational, and the normalization map is bijective.
\end{proof}

\medskip

\noindent
{\bf Remark 1.} In characteristic zero, $C$ is isomorphic to $\mP^1$
by \cite[2.2, 2.6]{Br93} . This does not generalize to arbitrary
characteristics: let indeed $K$ be a field of characteristic $p\ge 3$
and let $G=\SL(2)$, with its Borel subgroup $B$ of upper triangular
matrices. Let $M$ be the space of homogeneous polynomials of degree
$p+2$ in two variables $x,y$, where $G$ acts by linear
substitutions. The line $Ky^{p+2}\subset M$ is $B$--invariant; let 
$N=M/Ky^{p+2}$ be the quotient $B$--module, and let $[x^{p+1}y]$ be
the image of $x^{p+1}y$ in $\mP(N)$. Then 
$$
C=\overline{B\cdot [x^{p+1}y]}
$$ 
is a $B$--invariant curve in $\mP(N)$. Since
$$
(x+ty)^{p+1} y = x^{p+1} y + tx^p y^2 + t^p x y^{p+1} + t^{p+1} y^{p+2},
$$
the curve $C$ is singular at $[xy^{p+1}]$, the unique $B$--fixed point
in $\mP(N)$. On the other hand, $C$ is a $B$--invariant curve in
$$
X=G\times^B\mP(N),
$$
and the latter is a simple, nonsingular $G$--projective variety.

\medskip

Another useful observation is the following

\begin{lemma}\label{finite}
Any simple $G$--projective variety contains only finitely many
$B$--invariant curves.  
\end{lemma}

\begin{proof}
Let $X$ be a simple $G$--projective variety. By Lemma \ref{rational},
it suffices to show that $X^{R(P)}$ contains only finitely many
$B$--invariant curves, where $R(P)$ is the radical of a minimal
parabolic subgroup $P\supset B$. Then the quotient group $P/R(P)$ is
isomorphic to $\SL(2)$ or $\PGL(2)$; it acts on $X^{R(P)}$ with a
unique closed orbit (since $R_u(P)\subset U$, and $X^U$ consists of a
unique point). Thus, we may assume that $G=\SL(2)$. Then every
$H^0(\omega)$ is multiplicity--free as a $T$--module, so that the set
$\mP(H^0(\omega))^T$ is finite. Thus, $X^T$ is finite by Lemma
\ref{point}. This implies our statement, by Lemma \ref{rational}
again.
\end{proof}

We may now state our second main result.

\begin{theorem}\label{curves}

\begin{enumerate}

\item
With the notation of Theorem \ref{divisors}, the sink of each $D_i$ is
a unique point $x_i^-$, isolated in $X^T$. Moreover, $D_i$ is the
closure of $X^-(x_i^-)$, whereas $X^+(x_i^-)=B\cdot x_i^-$; the
$B$--invariant curve $C_i=\overline{B\cdot x_i^-}$ intersects $D_i$ 
at the unique point $x_i^-$, and intersects no other $D_j$. As a
consequence, $(D_i\cdot C_i)$ is a positive integer, and 
$(D_i\cdot C_j)=0$ for $j\ne i$. 

\item
The cone $NE(X)$ is generated by the classes of $C_1,\ldots,C_r$, and
these form a basis of the rational vector space $N_1(X)_{\mQ}$. 

\item
Any negative stratum having an irreducible component of codimension
one is actually irreducible, with closure some $D_i$. Likewise, any
positive stratum having an irreducible component of dimension one is 
open in some $C_i$.

\item
If, in addition, $X$ is nonsingular, then every $D_i$ intersects
transversally $C_i$ at $x_i^-$. In particular, $D_i$ is nonsingular at
$x_i^-$, and $(D_i\cdot C_i)=1$. Moreover, all $C_i$ are isomorphic to
$\mP^1$, and their classes form a basis of the group $N_1(X)$.

\end{enumerate}

\end{theorem}

\begin{proof}
(1) Let $Y_i$ be the sink of $D_i$ and let $x\in Y_i$. Then $X^+(x)$
is positive--dimensional (since $x\neq x^-$), and 
$X^+(x)\cap D_i=D_i^+(x)=\{x\}$. Since $D_i$ is a Cartier divisor,
it follows that $X^+(x)$ is one--dimensional. As a consequence,
$\overline{B\cdot x}$ is a $B$--invariant curve with $T$--fixed points
$x$ and $x^-$. Now Lemma \ref{finite} implies that
$Y_i=\{x\}$. Since $X^+(Y_i)$ is $B$-invariant, we must have
$X^+(Y_i)=B\cdot x$.

Let $Y$ be an irreducible component of $X^T$ through $x$. Then
$X^+(Y)\cap D_i = D_i^+(Y\cap D_i)$ contains $x$ as an isolated point,
so that $\dim_x X^+(Y)\le 1$. On the other hand, $X^+(Y)$ contains the
curve $B\cdot y$ for any $y\in Y$; it follows that $Y=\{x\}$. In other
words, $x$ is isolated in $X^T$; hence $X^+(x)=B\cdot x$.

To show that $\overline{X^-(x)}=D_i$, we choose an open affine
$T$--invariant neighborhood $X(x)\subset X$ of $x$, where $D_i$ is a
principal divisor associated to a regular function $f$, eigenvector of
$T$. Then $X(x)$ contains both $X^-(x)$ and $X^+(x)=B\cdot x$;
moreover, the weight of $f$ is negative on $\lambda$ (since $f$
restricts to a non--zero function on $B\cdot x$). It follows that $f$
vanishes at any point of $X^-(x)$, so that $X^-(x)=D_i^-(x)$.
 
Since $x^-\notin D_i$, the curve $C_i=\overline{B\cdot x}$ intersects
$D_i$ at $x$ only. We show that $C_i$ intersects no $D_j$ for $j\neq i$.
Otherwise, $x\in D_j$. Moreover, $x$ is not the sink of $D_j$ by the
preceding step, so that $D_j^+(x)\neq \{x\}$. But
$D_j^+(x)\subseteq X^+(x)=B\cdot x$, and hence 
$B\cdot x\subseteq D_j$. Thus, $x^-\in D_j$, a contradiction.  

(2) follows from (1) together with Theorem \ref{divisors}.

(3) Let $X^-(Y)$ be a negative stratum having an irreducible component 
$D$ of codimension one in $X$. Then $D=D_i$ for some index
$i$. Together with (1), it follows that $x_i^-\in D$, whence
$Y=\{x_i^-\}$ and $X^-(Y)=D_i^-$. 

Next let $X^+(Y)$ be a positive stratum having an irreducible
component $C$ of dimension one. Then $C$ is $B$--invariant, so that
$C=B\cdot x$ for some $x\in Y$. Let $Z$ be an irreducible component of
$Y$ containing $x$, then $C\subseteq B\cdot Z\subseteq X^+(Y)$, and
$B\cdot Z$ is irreducible. Hence $C=B\cdot Z$, that is, $Y=\{x\}$. Now
every irreducible component of $X^+(x)$ is $B$--invariant and contains
$x$; thus, $X^+(x) = B\cdot x = C$. By ($*$), we have 
$\dim X\le \dim_x X^+(x)+\dim_x X^-(x)$. Thus, $X^-(x)$ is a
divisor at $x$. Since $X^-(x)$ is disjoint from $X^-$, it follows that 
$x=x_i^-$ for some index $i$, whence $C=C_i$.

(4) Since $X$ is nonsingular and $x_i^-$ is an isolated fixed point,
$X^-(x_i^-)$ and $X^+(x_i^-)$ intersect transversally at $x_i^-$;
thus, the same holds for their closures $D_i$ and $C_i$. Together with
(1) and Theorem \ref{divisors}, it follows that $D_i$ restricts to a
globally generated Cartier divisor of degree $1$ on $C_i$. But the
normalization of $C_i$ is isomorphic to $\mP^1$, so that
$C_i\cong\mP^1$. Finally, by (2), any $\gamma\in N_1(X)$ decomposes as 
$\gamma=\sum_{i=1}^r c_i C_i$, and every $c_i=(D_i\cdot\gamma)$ is an
integer.
\end{proof}

\medskip

\noindent
{\bf Remark 2.} In particular, any stratum of codimension at most one
in a simple $G$--projective variety contains a unique $T$--fixed
point. But the whole fixed point set may be infinite; in fact, it may
have arbitrary irreducible components, as shown by the following
construction.

Let $Y\subseteq \mP^n$ be a projective variety. Let $G=\SL(3)$ with
standard opposite Borel subgroups $B, B^-$; let $\omega$ be the
highest root, i.e., the highest weight of the adjoint
representation. Then the zero--weight space $H^0(\omega)^T$ has
dimension $2$, so that $H^0(n\omega)^T$ has dimension $\ge n+1$. 
Thus, we may regard $Y$ as a subvariety of $\mP(H^0(n\omega)^T)$.
Now let $X=\overline{G\cdot Y}$ (closure in $\mP(H^0(n\omega))$. 
Clearly, $X$ is a simple $G$--projective variety, and $Y$ is an
irreducible component of $X^T$. In fact, both $X^+(Y)$ and $X^-(Y)$
have codimension $3$ in $X$, since $U\cdot Y$ (resp.~$U^-\cdot Y$) is
dense in $X^+(Y)$ (resp.~$X^-(Y)$).

\medskip

Finally, we compare the group $N_1(X)$ with the Chow group $A_1(X)$ of
one--cycles modulo rational equivalence.

\begin{corollary}\label{equivalences}
For any normal, simple $G$--projective variety $X$, the map 
$A_1(X)\rightarrow N_1(X)$ is an isomorphism over the rationals. If,
in addition, $X$ is nonsingular, then this map is an isomorphism.
\end{corollary}

\begin{proof}
We begin with the case where $X$ is nonsingular. Then it suffices to
show that the classes of $C_1,\ldots, C_r$ generate the group
$A_1(X)$. Decomposing $X$ into the positive strata $X^+(Y)$ and using
the long exact sequence for Chow groups \cite[Proposition 1.8]{Fu98},
this reduces to checking the vanishing of $A_1(X^+(Y))$ whenever 
$\dim X^+(Y)\ge 2$. 

Recall that $p^+:X^+(Y)\rightarrow Y$ is an affine bundle. By Lemma
\ref{sink} and \cite[Proposition 1.9]{Fu98}, it follows that
$A_1(X^+(Y))=0$ unless $\dim X^+(x)=1$ for all $x\in Y$. In the latter
case, since $B\cdot x\subseteq X^+(x)$ for all $x\in Y$, the fibers of
$p^+$ are $B$--invariant curves in $X^+(Y)$. By Lemma \ref{finite}, it
follows that $Y$ is a unique point, so that $X^+(Y)$ is
one--dimensional. 

For arbitrary $X$, recall that the group $A_1(X)$ is generated by
classes of $B$--invariant curves. Together with Theorem \ref{curves}
(3), this reduces to checking the following assertion: let $X^+(Y)$ be
a positive stratum and let $C\subseteq X^+(Y)$ be a $B$--invariant
curve which is not an irreducible component. Then the class of $C$ is
zero in $A_1(X^+(Y))_{\mQ}$.

Since $C$ is contained in an irreducible component of $X^+(Y)$ of
dimension $\ge 2$, we may find a closed irreducible $B$--stable
surface $S\subseteq X^+(Y)$ containing $C$. Let $Z=p^+(S)$, then
$Z\neq S$ (otherwise $S\subseteq X^T$, whence $S\subseteq X^B$, a
contradiction). For any $x\in Z$, the fiber at $x$ of
$p^+:S\rightarrow Z$ contains $B\cdot x$. By Lemma \ref{finite}, it
follows that $Z$ consists of a unique point $x$, the source of $S$. In
particular, $S$ is affine, and its algebra of regular functions is
negatively graded. Consider now the normalization $\tilde S$ of
$S$. The group $B$ acts on $\tilde S$ without fixed points, so that
$\tilde S$ is nonsingular. On the other hand, the algebra of regular
functions of $\tilde S$ is negatively graded as well; as a
consequence, $\tilde S$ is isomorphic to the affine plane. In
particular, $A_1(\tilde S)=0$, whence $A_1(S)_{\mQ}=0$. This proves
our assertion. 
\end{proof}

\section{The simple complete $\PGL(2)$--embeddings}

As an illustration of our results, we describe the cones of effective
one--cycles and nef divisors for all {\it simple complete embeddings
of $\PGL(2)$}, that is, for the normal complete $\PGL(2)$--varieties
containing an open orbit isomorphic to $\PGL(2)$ and a unique closed
orbit. Along the way, we obtain a realization of these embeddings as
$\PGL(2)$--orbit closures in a product of copies of the projective
line, suggested by work of Iozzi and Poritz \cite{IP99}. 

Note that there is a combinatorial classification of all embeddings of
$\PGL(2)$, presented e.g.~in \cite{MJ90}, as part of the Luna--Vust
theory of embeddings of homogeneous spaces; the geometric realization
of simple complete $\PGL(2)$--embeddings can be deduced from this
classification. The Chow rings of smooth complete $\SL(2)$--embeddings, 
and their cones of effective one--cycles, are described in \cite{MJ92}.

We introduce some notation. Let $G=PGL(2)$ and let $T$ be the image in
$G$ of the torus of diagonal matrices in $\SL(2)$. The map 
$t\mapsto \left(\begin{matrix} t&0\cr 0&t^{-1}\cr\end{matrix}\right)$ 
defines a one--parameter subgroup $\lambda$ of $T$, which generates
the group $\cX_*(T)$. The opposite Borel subgroups $B=G(\lambda)$,
$B^-=G(-\lambda)$ are the images in $G$ of the subgroups of upper
(resp.~lower) triangular matrices in $\SL(2)$. 

For the standard action of $G$ on $\mP^1$, the isotropy group of the
sink $\infty$ (resp.~the source $0$) is $B$ (resp.~$B^-$). We will
consider the diagonal action of $G$ on the product $(\mP^1)^r$ of
$r\ge 3$ copies of $\mP^1$. This is a (very) simple $G$--projective
variety, with closed orbit being the small diagonal, $\diag\mP^1$; the
sink is $\infty^r$, with $B^-$--invariant open affine neighborhood 
$(\mP^1\setminus\{0\})^r\cong \mA^r$ where $T$ acts by scalar
multiplication. The irreducible divisors of Theorem \ref{divisors}
are the 
$$
D_i=(\mP^1)^{i-1}\times\{0\}\times (\mP^1)^{r-i},
$$ 
and the curves of Theorem \ref{curves} are the lines
$$
C_i=\infty^{i-1}\times \mP^1\times \infty^{r-i},
$$ 
where $1\le i\le r$.

\begin{proposition}\label{PGL}
Let $p_1,\ldots,p_r\in\mP^1$ be pairwise distinct points and let
$$
X=X(p_1,\ldots,p_r)\subseteq (\mP^1)^r
$$ 
be the closure of the orbit $G\cdot(p_1,\ldots,p_r)$. Then the
following statements hold. 

\begin{enumerate}
\item
The irreducible components of the boundary $X\setminus G(p_1,\ldots,p_r)$
are the divisors
$$
\partial_i X =  \{(x^{i-1},y,x^{r-i})~\vert~ x,y\in\mP^1\} 
= G\cdot C_i~ ~(1\le i\le r).
$$
These are isomorphic to $\mP^1\times\mP^1$, and they intersect along
$\diag \mP^1$. 

\item
The open subset $X^-$ is isomorphic to $\mA^1\times\Sigma$, where
$\Sigma$ is the affine cone over the rational normal curve in
$\mP^{r-2}$. As a consequence, $X$ is a normal, Cohen--Macaulay,
simple $G$--projective variety; it is singular along $\diag\mP^1$, if
$r\ge 4$.

\item
The irreducible divisors constructed in Theorem \ref{divisors} are the
pull--backs of $0$ under the $r$ projections $X\rightarrow\mP^1$; they
are normal and Cohen--Macaulay. The curves constructed in Theorem
\ref{curves} are the lines $C_1,\ldots, C_r$.

\item
The boundary divisors $\partial_1 X,\ldots,\partial_r X$ form a basis
of the divisor class group of $X$ with rational coefficients; the
latter is isomorphic to the Picard group with rational coefficients. 

\item
A canonical divisor for $X$ is
$$
-K_X = \partial_1 X + \cdots + \partial_r X 
= \frac{2}{r-2}(D_1+\cdots+D_r).
$$
As a consequence, $-(r-2)K_X$ is very ample, so that $X$ is $\mQ$--Fano.

\end{enumerate}

Conversely, any simple complete embedding of $G$ is isomorphic to 
$X(p_1,\ldots,p_r)$, where $r\ge 3$ and $p_1,\ldots,p_r\in\mP^1$ are 
uniquely determined up to permutation and diagonal action of $G$.

\end{proposition}

\begin{proof}
(1) If $r= 3$, then $X=(\mP^1)^3$ and our assertions are evident; so
we may assume that $r\ge 4$. We may also assume that $p_1=\infty$ and
$p_2=0$; then $p_3,\ldots,p_r$ are distinct non--zero scalars.

The first projection $z_1:X\rightarrow \mP^1=G\cdot p_1\cong G/B$ is a
locally trivial fibration. It yields an isomorphism 
$X\cong G\times^B S$, where 
$S=\overline{B\cdot(p_1,\ldots,p_r)}\subset X$ is an irreducible 
projective $B$--invariant surface with a unique fixed point
$\infty^r$. Moreover, the map 
$$
U^-\times S^-\rightarrow X^-,~(g,s)\mapsto g\cdot s
$$ 
is an isomorphism, where $S^-=S\cap(\mP^1\setminus\{0\})^r$ is an
irreducible surface in $\mA^r$, invariant under scalar
multiplication. Since $S^-\subseteq\{\infty\}\times\mA^{r-1}$ (as
$p_1=\infty$), we will regard $S^-$ in $\mA^{r-1}$ via projection. 

Since
$\left(\begin{matrix}t&u\cr 0&t^{-1}\cr\end{matrix}\right)\cdot p_i=
t^2p_i+tu$ in $\mP^1\setminus\{\infty\}$, we see that $S^-$ is the
closure of the image of the rational map
$$
\mA^2-\rightarrow \mA^{r-1},~ (t,u)\mapsto
(\frac{1}{tu},\frac{1}{t^2p_3+tu},\ldots,\frac{1}{t^2p_r+tu}).
$$
For $2\le i\le r$, let $x_i = (t^2p_i+tu)^{-1}$. Equivalently, 
$tux_2=1$, and $tp_ix_i+u(x_i-x_2)=0$ for $3\le i\le r$. Hence
$S^-$ is defined (as a closed subset of $\mA^{r-1}$) by the vanishing
of all $2\times 2$ minors of the $2\times (r-2)$ matrix
$$
A=\left(
\begin{matrix}
p_3 x_3 & p_4 x_4 & \cdots & p_r x_r \cr
x_3-x_2& x_4-x_2 & \cdots & x_r-x_2\cr
\end{matrix}\right).
$$

It follows that the boundary $S^-\setminus B\cdot (p_1,\ldots,p_r)$
consists of the diagonal 
$$
(x_2=x_3=\cdots=x_r)=S^-\cap G\cdot C_1,
$$
together with the $r-1$ coordinate lines
$$
(x_2=\cdots=x_{i-1}=x_{i+1}=\cdots=x_r=0)=C_i^-.
$$
On the other hand, the open subset $X^-$ equals $U^-\cdot S^-$, and
intersects all $G$--orbits in $X$ (since $X^-$ intersects the unique
closed orbit). Thus, the boundary $X\setminus G\cdot (p_1,\ldots,p_r)$
is the union of the irreducible divisors $\partial_1 X = G\cdot C_1,
\ldots,\partial_r X = G\cdot C_r$.

(2) Clearly, $X$ is a simple $G$--projective variety. We claim that
its regular locus intersects all irreducible components of the
boundary. To see this, consider the projection
$(z_i,z_j,z_k):(\mP^1)^r\rightarrow (\mP^1)^3$ where 
$1\le i<j<k\le r$. Then  the restriction of $(z_i,z_j,z_k)$ to $X$ is
birational, and bijective over $G\cdot(C_j\setminus\{\infty^r\})$. Now
our claim follows from Zariski's main theorem.

By that claim, $X$ is nonsingular in codimension one. Next we
show that $X$ is Cohen--Macaulay, or, equivalently, that
$S^-$ is Cohen--Macaulay. For this, consider the closed subscheme
$\Sigma\subseteq \mA^{r-1}$ associated with the ideal generated by the
$2\times 2$ minors of the matrix $A$. Then $S^-$ is the support of
$\Sigma$. On the other hand, we may regard $\Sigma$ as a closed
subscheme of the space $M_{2,r-2}$ of $2\times (r-2)$ matrices, an
affine space of dimension $2r-4$. Then $\Sigma$ is the
scheme--theoretic intersection of the closed subscheme of matrices of
rank at most one (a Cohen--Macaulay variety of dimension $r-1$) with a
linear space of dimension $r-1$. Since $\dim \Sigma=2$, it follows
that $\Sigma$ is Cohen--Macaulay. Moreover, one easily checks that the
intersection $\Sigma\cap(x_2 = 0)\subseteq \mA^{r-1}$ is reduced, of
dimension one. Thus, $x_2$ is a non--zero--divisor on $\Sigma$, and
the latter is reduced. It follows that $\Sigma=S^-$.

Thus, $X^-\cong U^-\times S^-\cong\mA^1\times\Sigma$, where $\Sigma$ is
the affine cone over a curve $C\subset\mP^{r-2}$. Clearly, $C$ is
rational and not contained in any hyperplane; moreover, $C$ is
projectively normal, since $\Sigma$ is Cohen--Macaulay and nonsingular
in codimension one. Thus, $C$ is a rational normal curve.

(3) follows from the fact that $X$ contains $C_1,\ldots,C_r$, together
with Theorem \ref{curves}.

(4) For the divisor class group $A_2(X)$, we have an exact sequence
$$
\bigoplus_{i=1}^r \mZ \; \partial_i X \rightarrow A_2(X)\rightarrow 
A_2(G\cdot(p_1,\ldots,p_r))\rightarrow 0.
$$
Moreover, the group $A_2(G\cdot(p_1,\ldots,p_r))=\Pic(\PGL(2))$ has
order $2$, and the classes of the $\partial_i X$ are linearly
independent in $A_2(X)_{\mQ}$ (since every regular invertible function
on $\PGL(2)$ is constant). Therefore, these classes form a basis of
$A_2(X)_{\mQ}$. The lattter contains $\Pic(X)_{\mQ}$ as a subspace of
dimension $r$; thus, they are equal.

(5) By (4), we have 
$$
-K_X=a_1\partial_1 X + \cdots + a_r\partial_r X
$$ 
for unique rational coefficients $a_1,\ldots,a_r$. If $r=3$, then
$X=(\mP^1)^3$ and one obtains easily $a_1=a_2=a_3=1$. In the general
case, any $\partial_i X$ is mapped isomorphically to its image under
some projection to $(\mP^1)^3$; it follows that $a_i=1$.

To express $-K_X$ in terms of $D_1,\ldots,D_r$, consider the divisors
of the rational functions $z_i - z_j$, where
$z_1,\ldots,z_r:X\rightarrow \mP^1$ denote the projections. We obtain
$$
\div(z_i-z_j)=-D_i - D_j + \sum_{k, k\notin \{i, j\}}\partial_k X,
$$
whence 
$$
\partial_1 X + \cdots + \partial_r X =\frac{2}{r-2}(D_1+\cdots+D_r)
$$
in $A_2(X)_{\mQ}$.

For the final assertion, let $X$ be a simple complete embedding of $G$. 
Choose $p\in X$ in the open $G$--orbit. By Lemma \ref{point}, there
exists a finite equivariant morphism
$$
f:X\rightarrow\mP(H^0(\omega)).
$$ 
Now $H^0(\omega)$ is the space of homogeneous polynomials of degree
$d$ in two variables $x,y$, where $G$ acts by linear substitutions;
here $d$ is an even positive integer (see
e.g.~\cite[II.2.16]{Ja87}). Let $P(x,y)\in H^0(\omega)$ be a
representative of $f(p)\in\mP(H^0(\omega))$; let
$$
P(x,y)=\prod_{i=1}^s (a_i x+b_i y)^{m_i}
$$ 
be a decomposition into a product of pairwise distinct linear forms,
with multiplicities. Since $f$ is finite, then so is the isotropy
group of the line $KP$. Thus, $s\ge 3$.

The source of the $\lambda$--action on $\mP(H^0(\omega))$ consists of
the image of $y^d$, and the corresponding open subset
$\mP(H^0(\omega))^-$ consists of the images of those homogeneous
polynomials where the coefficient of $y^d$ is non--zero. Since $f$ is
finite and equivariant, $X^-=f^{-1}(\mP(H^0(\omega))^-$. It follows
that $(X\setminus X^-)\cap G\cdot p$ consists of $s$ irreducible
components. On the other hand, any irreducible component $D_i$ of
$X\setminus X^-$ intersects $G\cdot p$ (otherwise, $D_i$
is $G$--invariant, and hence contains the closed $G$--orbit, a
contradiction). Thus, $s=r$ with the notation of Theorem
\ref{divisors}.

Next let 
$$
f_i:X\rightarrow Y_i=\Proj\bigoplus_{n=0}^{\infty} H^0(X,\cO_X(nD_i))
$$ 
be the morphism associated with the globally generated divisor $D_i$.
Then $Y_i$ is a normal simple $G$--projective variety, and $f_i$ is
equivariant, separable, with connected fibers. Moreover, the Picard
group of $Y_i$ is freely generated by a unique ample divisor, which
pulls back to $D_i$. By considering a finite morphism 
$Y_i\rightarrow \mP(H^0(\omega_i))$ and arguing as above, one obtains
that $\dim Y_i =1$. On the other hand, $G$ acts on $Y_i$ with a dense
separable orbit and a unique closed orbit, whence $Y_i$ is isomorphic
to the projective line with standard $G$--action. 

Consider now the product morphism 
$$
f=\prod_{i=1}^r f_i:X\rightarrow\prod_{i=1}^r Y_i \cong (\mP^1)^r.
$$
Then $f$ is finite, since $D_1+\cdots+D_r$ is ample. Let
$f(p)=(p_1,\ldots,p_r)$, then the $p_i$ are pairwise distinct (since
so are the $D_i$), and their number is at least $3$. Thus, the
restriction $G\cdot p\rightarrow G\cdot f(p)$ is an isomorphism. On
the other hand, $f(X)$ is normal by (1), so that $f$ is a closed
immersion by Zariski's main theorem.

So we have proved that $X$ embeds equivariantly into $(\mP^1)^r$,
where $r$ is the rank of $\Pic(X)$; moreover, the $r$ projections
$z_i:X\rightarrow\mP^1$ are uniquely determined. This implies the
remaining uniqueness assertion.
\end{proof}

\end{document}